# Market Clearing for Uncertainty, Generation Reserve, and Transmission Reserve—Part II: Case Study

Hongxing Ye, *Student Member, IEEE,* Yinyin Ge, *Student Member, IEEE,*
Mohammad Shahidehpour, *Fellow, IEEE,* Zuyi Li, *Senior Member, IEEE*

*Abstract*—In Part II of this two-part paper, we analyze the marginal prices derived in Part I of this two-part paper within a robust optimization framework. The load and generation are priced at Locational Marginal Price (LMP) while the uncertainty and generation reserve are priced at Uncertainty Marginal Price (UMP). The Financial Transmission Right (FTR) underfunding is demonstrated when there is transmission reserve. A comparison between traditional reserve price and UMP is presented. We also discuss the incentives for market participants within the new market scheme.

## I. Introduction

IN Part I of this two-part paper, the theory for market clearing is presented [1]. The case studies are reported in Part II of this two-part paper. Two systems are simulated in this paper. The first one is a 6-bus system, and the second one is the IEEE 118-bus system. For the 6-bus system, the basic ideas of UMP are presented within the robust optimization framework. We show how to use the UMPs in the new scheme to clear the market in Section II-A. A comparison between the UMP and traditional reserve price is made in Section II-B without transmission constraint. By dropping the most challenging transmission constraint, the robustness of the traditional reserve scheme can also be guaranteed. Thus, we can compare the new scheme and the traditional scheme fairly. In Section II-C, the sensitivity analysis regarding UMP is performed. We study the impacts of different uncertainty levels, ramping rates, generation capacities, and energy bids on UMPs. In Section II-D, the FTR underfunding issue is presented within the robust optimization framework. For the IEEE 118-Bus system, the UMPs related products are presented for different uncertainty levels in Section III-A. We also analyze the behaviors and impacts of flexible sources using an energy storage example in Section III-B.

## II. 6-bus System

A 6-bus system is studied in this section. The one-line diagram is shown in Fig.1. The unit data and line data are shown in Table I and Table II, respectively. Table III

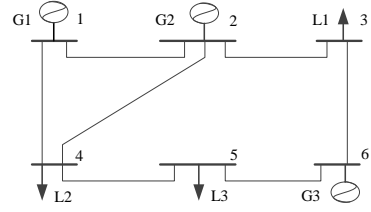

Fig. 1. One-line Diagram for 6-bus system

TABLE I
Unit Data for the 6-bus System

| # | $P^{\min}$ | $P^{\max}$ | $P_0$ | a | b | c | $R^u$ | $R^d$ | $C_u$ | $C_d$ | $T^{\text{on}}$ | $T^{\text{off}}$ | $T_0$ |
|---|---|---|---|---|---|---|---|---|---|---|---|---|---|
| 1 | 100 | 220 | 120 | 0.004 | 13.5 | 176.9 | 24 | 24 | 180 | 50 | 4 | 4 | 4 |
| 2 | 10 | 100 | 50 | 0.001 | 32.6 | 129.9 | 12 | 12 | 360 | 40 | 3 | 2 | 3 |
| 6 | 10 | 20 | 0 | 0.005 | 17.6 | 137.4 | 5 | 5 | 60 | 0 | 1 | 1 | −2 |

$P^{\min}, P^{\max}, P_0$: min/max/initial generation level (MW); fuel cost ($): $aP^2 + bP + c$; $R^u, R^d$: ramping up/down rate (MW/h); $C_u, C_d$: startup/shutdown cost ($); $T^{\text{on}}, T^{\text{off}}, T_0$: min on/min off/initial time (h) presents the load and uncertainty information. Column "Base Load" shows the hourly forecasted load. Assume that the load distributions are 20%, 40%, and 40% for Bus 3, Bus 4, and Bus 5, respectively. $\bar{u}_{1,t}$ and $\bar{u}_{3,t}$ in Table III are the bounds of the uncertainties at Bus 1 and Bus 3, respectively. The uncertainty bounds at other buses are 0. It is assumed that the relative forecasting errors increase with hours. Uncertainty $\epsilon_{1,t}$ and $\epsilon_{3,t}$ also respect

$$-\Lambda \cdot \bar{u}_{m,t} \le \epsilon_{m,t} \le \Lambda \cdot \bar{u}_{m,t}, \forall t, m \tag{1a}$$

$$\sum_m \frac{|\epsilon_{m,t}|}{\bar{u}_{m,t}} \le \Lambda_\Delta, \forall t, \tag{1b}$$

where (1a) denotes the uncertainty interval at a single bus, and (1b) represents the system-wide uncertainty [2], [3]. The $\Lambda$ and $\Lambda^\Delta$ are the budget parameters for the single bus and system, respectively.

TABLE II
Line Data for the 6-bus System

| from | 1 | 1 | 2 | 5 | 3 | 2 | 4 |
|---|---|---|---|---|---|---|---|
| to | 2 | 4 | 4 | 6 | 6 | 3 | 5 |
| x(p.u.) | 0.17 | 0.258 | 0.197 | 0.14 | 0.018 | 0.037 | 0.037 |
| capacity(MW) | 200 | 100 | 100 | 100 | 100 | 200 | 200 |

TABLE III
LOAD AND UNCERTAINTY DATA FOR THE 6-BUS SYSTEM (MW)

| Time (h) | Base Load | $\bar{u}_{1,t}$ | $\bar{u}_{3,t}$ | Time (h) | Base Load | $\bar{u}_{1,t}$ | $\bar{u}_{3,t}$ |
|---|---|---|---|---|---|---|---|
| 1 | 175.19 | 1.09 | 0.29 | 13 | 242.18 | 19.68 | 5.25 |
| 2 | 165.15 | 2.06 | 0.55 | 14 | 243.6 | 21.32 | 5.68 |
| 3 | 158.67 | 2.98 | 0.79 | 15 | 248.86 | 23.33 | 6.22 |
| 4 | 154.73 | 3.87 | 1.03 | 16 | 255.79 | 25.58 | 6.82 |
| 5 | 155.06 | 4.85 | 1.29 | 17 | 256 | 27.2 | 7.25 |
| 6 | 160.48 | 6.02 | 1.6 | 18 | 246.74 | 27.76 | 7.4 |
| 7 | 173.39 | 7.59 | 2.02 | 19 | 245.97 | 29.21 | 7.79 |
| 8 | 177.6 | 8.88 | 2.37 | 20 | 237.35 | 29.67 | 7.91 |
| 9 | 186.81 | 10.51 | 2.8 | 21 | 237.31 | 31.15 | 8.31 |
| 10 | 206.96 | 12.94 | 3.45 | 22 | 232.67 | 31.99 | 8.53 |
| 11 | 228.61 | 15.72 | 4.19 | 23 | 195.93 | 28.16 | 7.51 |
| 12 | 236.1 | 17.71 | 4.72 | 24 | 195.6 | 29.34 | 7.82 |

TABLE IV
MARGINAL COSTS AT DIFFERENT GENERATION LEVELS ($/MWH)

| Gen. 1 | | | Gen. 2 | | | Gen. 3 | | |
|---|---|---|---|---|---|---|---|---|
| $\underline{P}_1^w$ | $\bar{P}_1^w$ | mar. cost | $\underline{P}_2^w$ | $\bar{P}_2^w$ | mar. cost | $\underline{P}_3^w$ | $\bar{P}_3^w$ | mar. cost |
| 100 | 124 | 14.396 | 10 | 28 | 32.638 | 10 | 12 | 17.71 |
| 124 | 148 | 14.588 | 28 | 46 | 32.674 | 12 | 14 | 17.73 |
| 148 | 172 | 14.78 | 46 | 64 | 32.71 | 14 | 16 | 17.75 |
| 172 | 196 | 14.972 | 64 | 82 | 32.746 | 16 | 18 | 17.77 |
| 196 | 220 | 15.164 | 82 | 100 | 32.782 | 18 | 20 | 17.79 |

### A. LMP and UMP

In this subsection, we show how to obtain the marginal prices for energy, uncertainty, and generation reserve according to Eqs. (11) and (14) in [1]. The payments and credits based on these prices are also calculated for market participants.

Consider the case where $\Lambda = 1, \Lambda^\Delta = 2$. The CG based approach converges after 2 iterations. Hence, $\mathcal{K} = \{1, 2\}$. The UC results are shown in Fig. 2. The ON/OFF statuses of G1, G2, and G3 are denoted by black, blue, and red dots, respectively. It can be observed that the cheaper unit G1 is always on, G2 is off from Hour 5 to Hour 10, and G3 is off from Hour 1 to Hour 9. Given the UC solutions, the problem (RSCED) in [1] can be solved by commercial linear programming (LP) solver. The marginal prices are then obtained as byproducts.

The generation outputs are presented in Table V at Hours 21 and 22. It can be observed that G1 supplies most of the loads at Hour 21, which is 195.19 MW. According to the bid information in Table IV, G2 is much more expensive than G1 and G3. Hence, the output of G2 is relatively small and at the low level of its capacity. The upward and downward generation reserves provided by the three units are also listed in Table V. These data can be obtained directly from Eqs. (18) and (19) given the generation output $P_{i,t}$. The maximum upward reserve is 24MW provided by G1, and the minimum one is 3.46 MW provided by G3. Although the remaining generation capacity of G1 is $220 − 195.19 = 24.62$ MW, the upward reserve is limited by its upward ramping rate 24 MW. In the meantime, the upward reserve provided by G3 is limited by its generation capacity although it has more remaining ramping capacity (i.e. $\min\{20 − 16.54, 5\} = 3.46$ MW).

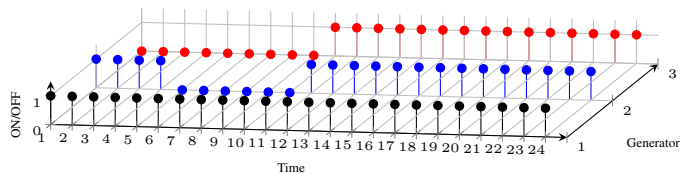

Fig. 2. Unit Commitment Results ($\Lambda = 1, \Lambda^\Delta = 2$)

TABLE V
GENERATION OUTPUT AND RESERVE ($\Lambda = 1, \Lambda^\Delta = 2$, MW)

| T | $P_1$ | $P_2$ | $P_3$ | $Q_1^{\text{up}}$ | $Q_1^{\text{down}}$ | $Q_2^{\text{up}}$ | $Q_2^{\text{down}}$ | $Q_3^{\text{up}}$ | $Q_3^{\text{down}}$ |
|---|---|---|---|---|---|---|---|---|---|
| 21 | 195.19 | 25.58 | 16.54 | 24 | −24 | 12 | −12 | 3.46 | −5 |
| 22 | 188.26 | 28.93 | 15.48 | 24 | −24 | 12 | −12 | 4.52 | −5 |

TABLE VI
EXTREME POINTS OF UNCERTAINTY SET

| | $k=1$ | | | | $k=2$ | | | |
|---|---|---|---|---|---|---|---|---|
| $t$ | $\hat{\epsilon}^1_{1,t}$ | $\hat{\epsilon}^1_{3,t}$ | $\pi^{u,1}_{1,t}$ | $\pi^{u,1}_{3,t}$ | $\hat{\epsilon}^2_{1,t}$ | $\hat{\epsilon}^2_{3,t}$ | $\pi^{u,2}_{1,t}$ | $\pi^{u,2}_{3,t}$ |
| 21 | 31.15 | 8.31 | 14.87 | 14.87 | −31.15 | 8.31 | −17.67 | 1.77 |
| 22 | −31.99 | −8.53 | 0 | 0 | 31.99 | 8.53 | 29.81 | 29.81 |

Table VI shows the extreme points obtained in the CG-based approach. The intermediate price signals for these points $\pi^{u,k}_{m,t}$ are also presented. It can be observed that the worst point is always obtained at the extreme point of the uncertainty set. For example, at Hour 21, the $\hat{\epsilon}^1_{1,1}$ is 31.15 MW. It is exactly the upper bound of the uncertainty at Hour 21 at Bus 1. The data in Table VI also verifies Lemma 1. The intermediate UMPs $\pi^{u,k}_{m,t}$ have the same sign as the uncertainties $\hat{\epsilon}^k_{m,t}$ at the same bus. At $t = 22, k = 1$, the uncertainties are negative while $\pi^{u,1}_{1,22} = \pi^{u,1}_{3,22} = 0$.

The LMPs, aggregated upward UMPs, and aggregated downward UMPs are shown in Table VII, Table VIII, and Table IX, respectively. It is noted that UMPs still exist at buses without uncertainties (i.e., Buses 2,4,5,6). This is similar to LMPs, which also exist at buses where net power injections are 0. As shown in Table VII, the LMPs vary with locations at Hour 21. It indicates that the line congestion exists at Hour 21. The load at Bus 4 has to pay the highest LMP $43.71/MWh. The UMPs at Hour 21 are also different at various locations. The highest upward UMP at Hour 21 is also located at Bus 4. As shown in Table IX, the downward UMPs are zeros except the one at Bus 1 at Hour 21. With these prices, the market participants can be paid and credited. For example, G1 supplies 195.19 MWh load at $14.97/MWh at Hour 21. Then the energy payment to G1 is

$$195.19\text{MWh} \times \$14.97/\text{MWh} = \$2922.99.$$

At the same time, G1 also provides the 24 MW upward generation reserve and the -24 MW downward generation reserve. Hence, the generation reserve credit entitled to G1 is

$$24\text{MW} \times \$14.87/\text{MW} + (-24)\text{MW} \times \$(-17.67)/\text{MW} = \$780.96,$$

TABLE VII
LOCATIONAL MARGINAL PRICE ($\Lambda = 1, \Lambda^{\Delta} = 2$,$/MWH)

| t | Bus1 | Bus2 | Bus3 | Bus4 | Bus5 | Bus6 |
|---|------|------|------|------|------|------|
| 21 | 14.97 | 32.64 | 34.4 | 43.71 | 41.94 | 35.26 |
| 22 | 47.56 | 47.56 | 47.56 | 47.56 | 47.56 | 47.56 |

TABLE VIII
UPWARD UMP ($\Lambda = 1, \Lambda^{\Delta} = 2$,$/MW)

| t | Bus1 | Bus2 | Bus3 | Bus4 | Bus5 | Bus6 |
|---|------|------|------|------|------|------|
| 21 | 14.87 | 14.87 | 16.63 | 25.94 | 24.17 | 17.49 |
| 22 | 29.81 | 29.81 | 29.81 | 29.81 | 29.81 | 29.81 |

TABLE IX
DOWNWARD UMP ($\Lambda = 1, \Lambda^{\Delta} = 2$,$/MW)

| t | Bus1 | Bus2 | Bus3 | Bus4 | Bus5 | Bus6 |
|---|------|------|------|------|------|------|
| 21 | $-17.67$ | 0 | 0 | 0 | 0 | 0 |
| 22 | 0 | 0 | 0 | 0 | 0 | 0 |

where the negative sign of UMP means it is downward price, and the downward generation reserve is also negative. Hence, the credit is positive. The load located at Bus 3 is charged

$$47.464 \text{MWh} \times \$34.4/\text{MWh} = \$1632$$

for energy. The uncertainty source at Bus 3 is charged

$$8.31 \text{MW} \times \$16.63/\text{MW} = \$138.19.$$

The LMPs and UMPs also vary with the time intervals. It is observed that the LMPs at different buses are the same at Hour 22, which means that there is no line congestion. The downward UMPs at Hour 22 are zeros for all buses, and the upward UMPs are $29.81/MW.

Next, we analyze if the generators are inclined to deviate from the dispatch instructions given the LMP and UMP signals. At Hour 21, the LMPs are the same as the marginal costs of G1 and G2 according to the data in Table VII and Table IV. It means that both G1 and G2 cannot get more profits on energy by increasing or decreasing the outputs. The reserves provided by G1 and G2 are limited by their physical ramping rates. Hence, they can do nothing to increase the generation reserve payments. On the other hand, the marginal cost for G3 is $17.77/MWh when the output is $16.54 \in [16, 18]$. The LMP paid to G3 is $35.26/MWh on Bus 6, which is $17.49/MWh larger than its marginal cost. It seems that if G3 provides more power, then it would get more profits on energy. However, increasing its output also results in lower upward generation reserve according to Eq. (18) in [1]. The upward UMP is $17.49/MW on Bus 6, which is exactly the difference between the LMP and G3's marginal cost. In other words, G3 cannot get more profits by deviating from the dispatch instruction. Hence, the UMP on Bus 6 is also the opportunity cost for G3 [4]. The analysis here also verifies the competitive market equilibrium shown in Section III-D of the companion paper.

In fact, the UMPs provide important price signals on the planning of renewable energy sources and storages. For example, according to the above analysis, the UMP at Bus 2 is relatively small. It is an ideal location for renewable energy sources in terms of payment for uncertainties. On the other hand, the UMP at Bus 4 is large, which may attract the long-term investment for storages or generation plants with large ramping rates because they can get more profits by providing generation reserves. As the UMPs are derived from the Lagrangian function, these planning activities also lower the total operation cost and increase the social welfare.

### B. Comparison between Existing LMPs and Reserve Prices

The motivation of this subsection is to compare the proposed clearing scheme with the existing one. However, as the reserve is not robust in the traditional scheme, we cannot compare them fairly. With the observation that the transmission constraints are the most challenging one in the robust SCUC framework, we drop these constraints in this subsection. Doing so is equivalent to assuming that transmission line limits are large enough so that no line flow constraint is binding. By relaxing the transmission constraints, the robustness of the tranditional UC solution can be guranteed by adding the following spinning reserve constraints

$$I_{i,t} P_i^{\min} \leq Q_{i,t}^{\text{down}} + P_{i,t}, \quad Q_{i,t}^{\text{up}} + P_{i,t} \leq I_{i,t} P_i^{\max}, \forall i, t \quad (2a)$$

$$-R_i^d I_{i,t} \Delta T \leq Q_{i,t}^{\text{down}}, \quad Q_{i,t}^{\text{up}} \leq R_i^u I_{i,t} \Delta T, \forall i, t \quad (2b)$$

$$\sum_i Q_{i,t}^{\text{down}} \leq \underline{R}_t, \quad \sum_i Q_{i,t}^{\text{up}} \geq \bar{R}_t, \forall t, \quad (2c)$$

where $Q_{i,t}^{\text{up}}$ and $Q_{i,t}^{\text{down}}$ is the largest upward and downward reserve, respectively. $\underline{R}_t$ and $\bar{R}_t$ are system wide reserve requirement. Refer to [5], [6] for more details on the reserve formulations. In this paper, $\Delta T$ is set to 1. The reserve requirements $\underline{R}_t$ and $\bar{R}_t$ are set to the system-wide lower bound and upper bound uncertainty in (1b), respectively. By adding the above (2a-2c) to (1a-3c) in [1], we can form a new problem, whose solution is immunized against all the uncertainties and the same as the solution to RSCUC without transmission constraints. Therefore, we are able to compare the proposed pricing mechanism with the traditional one.

Consider the case where $\Lambda = 0.8, \Lambda^{\Delta} = 2$. The generation schedules obtained from RSCUC and traditional SCUC are presented in Table X. The optimal solutions to the two models are the same. It verifies the analysis shown in the first paragraph in this subsection. It indicates that when the upper bound and lower bound of the system-wide uncertainty are available, the UC and generation schedules are robust against uncertainties in the traditional SCUC by setting the reserve requirement. Consider Hour 21, when the upper bound of the system-wide uncertainty is

$$(31.15 + 8.31) \times 0.8 = 31.57 \text{MW}.$$

At the same Hour 21, the upward generation reserve provided by the three units are

$$14.57 + 12 + 5 = 31.57 \text{MW}.$$

Therefore, the uncertainties can be accommodated by these reserves in the upward direction. Similar conclusion can be verified in the downward direction.

TABLE X
GENERATION OUTPUT AND RESERVE W/O T.C. (MW)

| T | $P_1$ | $P_2$ | $P_3$ | $Q_1^{\text{up}}$ | $Q_1^{\text{down}}$ | $Q_2^{\text{up}}$ | $Q_2^{\text{down}}$ | $Q_3^{\text{up}}$ | $Q_3^{\text{down}}$ |
|---|---|---|---|---|---|---|---|---|---|
| 21 | 205.432 | 16.878 | 15 | 14.568 | $-24$ | 12 | $-6.878$ | 5 | $-5$ |
| 22 | 204.254 | 18.416 | 10 | 15.746 | $-24$ | 12 | $-8.416$ | 5 | 0 |

TABLE XI
MARGINAL PRICES AND OPPORTUNITY COSTS W/O T.C.

| t | LMP | UMP/Res.P | $\pi_{1,t}^{o,\text{up}}$ | $\pi_{1,t}^{o,\text{down}}$ | $\pi_{2,t}^{o,\text{up}}$ | $\pi_{2,t}^{o,\text{down}}$ | $\pi_{3,t}^{o,\text{up}}$ | $\pi_{3,t}^{o,\text{down}}$ |
|---|---|---|---|---|---|---|---|---|
| 21 | 32.638 | 17.474 | 17.474 | 0 | 0 | 0 | 0 | 0 |
| 22 | 15.164 | $-17.474$ | 0 | 0 | 0 | $-17.474$ | 0 | $-17.474$ |

Table XI presents the LMPs and UMPs in the case without transmission constraints. The LMPs obtained from the new scheme and the traditional one are the same. The UMPs are also the same as the reserve prices obtained from the traditional approach, which are the Lagrangian multipliers for (2c) in the simplified model [6], [7]. As the transmission constraints are dropped, the marginal prices do not change with the locations, and there is a unique system-wide LMP and UMP at each interval. It is noted that the downward UMP at Hour 21 is zero, and the upward UMP at Hour 22 is zero. In Table XI, only non-zero UMPs are listed. Consider the payment and credit related to uncertainties at Hours 21 and 22. The three units get the same payment based on the two schemes,

G1: $14.57 \times 17.47 + (-24 \times -17.47) = \$673.94$

G2: $12 \times 17.47 + (-8.42 \times -17.47) = \$356.75$

G3: $5 \times 17.47 + 0 = \$87.37$.

In the clearing scheme proposed in this paper, the uncertainty source at Bus 1 pays

$17.47 \times 31.15 \times 0.8 + (-17.47) \times (-31.99) \times 0.8 = \$882.65$

The uncertainty source at Bus 3 pays

$17.47 \times 8.31 \times 0.8 + (-17.47) \times (-8.53) \times 0.8 = \$235.36$

It can be observed that the total reserve credit $\$1118 = 673.94 + 356.75 + 87.37$ actually originates from the uncertainty payment $\$1118 = 882.65 + 235.36$. However, in the traditional scheme, the reserve cost is allocated to the load serving entity [7]. In means that L2 and L3, which are located at buses 4 and 5, respectively, pay $\$447.24 = 1118.1 \times 0.4$ less at Hour 21 in the new scheme than in the traditional one.

Another observation is that the UMPs in Table XI are determined by the largest opportunity costs of generators. In this case, the upward and downward opportunity costs for generators are calculated as

$$\pi_{i,t}^{o,\text{up}} := \sum_{k \in \mathcal{K}_{m,t}^{\text{up}}} \bar{\beta}_{i,t}^k, \quad \pi_{i,t}^{o,\text{down}} := \sum_{k \in \mathcal{K}_{m,t}^{\text{down}}} \underline{\beta}_{i,t}^k.$$

At Hour 21, the opportunity costs of G2 and G3 are zeros. It can also be verified by comparing the LMP and the marginal fuel costs. For example, generation schedule of G2 is 16.878 MW at Hour 21, and the LMP is \$32.638/MWh, which is the same as the marginal cost. In the meantime, the opportunity cost of G1 is \$17.474/MW, which is also the UMP. At Hour 22, the downward UMP is \$17.474/MW, which is the opportunity cost of G2 or G3.

When considering transmission constraints, the generation reserve cannot be guaranteed at bus levels in the traditional SCUC model. The marginal prices are different in the two schemes.

*C. Sensitivity Analysis*

In this subsection, we perform the sensitivity analysis of UMPs with respect to the uncertainty levels, unit capacities and ramping rates, and energy bids. The bus-level uncertainty level is defined by $\Lambda$, and the system uncertainty level is defined by the budget parameter $\Lambda^\Delta$. We have $\Lambda \in \{0.5, 0.8, 1\}$ and $\Lambda^\Delta \in \{1, 2\}$.

The total operation costs are presented in Table XII. It can be observed that when $\Lambda^\Delta$ is fixed, the cost is a monotonically non-decreasing function of $\Lambda$. For example, the cost is \$88,663 if $\Lambda^\Delta = 2$ and $\Lambda = 0.8$. The cost is increased to \$89,851 if $\Lambda$ is changed to 1. When $\Lambda$ is fixed, the cost is also a monotonically non-decreasing function of $\Lambda^\Delta$. The cost in the case of $\Lambda^\Delta = 2, \Lambda = 0.8$ is \$663 higher than that in the case of $\Lambda^\Delta = 1, \Lambda = 0.8$. We can draw the conclusion that the operation cost for the base case may increase with the uncertainty level. A special case is simulated with $\Lambda^\Delta = 0, \Lambda = 0$. It can be observed that the total operation cost is \$87,975, which is the same as that in the case of $\Lambda^\Delta = 1, \Lambda = 0.5$ and $\Lambda^\Delta = 2, \Lambda = 0.5$. It indicates that even if we ignore the uncertainties in the SCUC and SCED problems, certain level of robustness is still kept in the optimal solution as generation reserves are already available at individual buses. If the uncertainty level is below $\Lambda = 0.5$ at each bus, no additional cost is needed to immunize the system against the uncertainties. LMPs with respect to different uncertainty levels are presented in Table XIII. It can be observed that the LMPs also vary with different uncertainties. For example, all the LMPs from Bus 3 to Bus 6 increase with the increment of $\Lambda$, but LMP at Bus 1 decreases. It indicates that the uncertainty level is also reflected in the LMPs.

TABLE XII
OPERATION COST W.R.T. UNCERTAINTY LEVEL (\$)

| $\Lambda^\Delta$ | $\Lambda$ | Cost |
|---|---|---|
| 2 | 1 | 89,851 |
| 2 | 0.8 | 88,663 |
| 2 | 0.5 | 87,975 |
| 1 | 1 | 89,196 |
| 1 | 0.8 | 88,000 |
| 1 | 0.5 | 87,975 |
| 0 | 0 | 87,975 |

TABLE XIII
LMP W.R.T. UNCERTAINTY LEVEL ($t = 21$, \$/MWH)

| $\Lambda^\Delta$ | $\Lambda$ | Bus1 | Bus2 | Bus3 | Bus4 | Bus5 | Bus6 |
|---|---|---|---|---|---|---|---|
| 2 | 1 | 14.972 | 32.638 | 34.404 | 43.709 | 41.943 | 35.263 |
| 2 | 0.8 | 15.164 | 32.638 | 34.384 | 43.589 | 41.842 | 35.234 |
| 2 | 0.5 | 15.164 | 32.638 | 34.384 | 43.589 | 41.842 | 35.234 |
| 1 | 1 | 15.164 | 32.638 | 34.384 | 43.589 | 41.842 | 35.234 |
| 1 | 0.8 | 15.164 | 32.638 | 34.384 | 43.589 | 41.842 | 35.234 |
| 1 | 0.5 | 15.164 | 32.638 | 34.384 | 43.589 | 41.842 | 35.234 |

TABLE XIV
GENERATION OUTPUT AND RESERVE W.R.T UNCERTAINTY LEVEL ($t = 21$, MW)

| $\Lambda^\Delta$ | $\Lambda$ | $P_1$ | $P_2$ | $P_3$ | $Q_1^{\text{up}}$ | $Q_2^{\text{up}}$ | $Q_3^{\text{up}}$ | $Q_1^{\text{down}}$ | $Q_2^{\text{down}}$ | $Q_3^{\text{down}}$ |
|---|---|---|---|---|---|---|---|---|---|---|
| 2 | 1   | 195.193 | 25.577 | 16.54 | 24     | 12 | 3.46 | $-24$ | $-12$    | $-5$ |
| 2 | 0.8 | 201.589 | 20.721 | 15    | 18.411 | 12 | 5    | $-24$ | $-10.721$ | $-5$ |
| 2 | 0.5 | 203.173 | 14.137 | 20    | 16.827 | 12 | 0    | $-24$ | $-4.137$  | $-5$ |
| 1 | 1   | 196.023 | 21.287 | 20    | 23.977 | 12 | 0    | $-24$ | $-11.287$ | $-5$ |
| 1 | 0.8 | 202.253 | 15.057 | 20    | 17.747 | 12 | 0    | $-24$ | $-5.057$  | $-5$ |
| 1 | 0.5 | 203.173 | 14.137 | 20    | 16.827 | 12 | 0    | $-24$ | $-4.137$  | $-5$ |

TABLE XV
UPWARD UMP W.R.T. UNCERTAINTY LEVEL ($t = 21$, $/MW)

| $\Lambda^\Delta$ | $\Lambda$ | Bus1 | Bus2 | Bus3 | Bus4 | Bus5 | Bus6 |
|---|---|---|---|---|---|---|---|
| 2 | 1   | 14.868 | 14.868 | 16.634 | 25.939 | 24.173 | 17.493 |
| 2 | 0.8 | 0 | 0 | 1.746 | 10.951 | 9.204 | 2.596 |
| 2 | 0.5 | 0 | 0 | 0 | 0 | 0 | 0 |
| 1 | 1   | 0 | 0 | 1.746 | 10.951 | 9.204 | 2.596 |
| 1 | 0.8 | 0 | 0 | 1.746 | 10.951 | 9.204 | 2.596 |
| 1 | 0.5 | 0 | 0 | 0 | 0 | 0 | 0 |

TABLE XVI
DOWNWARD UMP W.R.T. UNCERTAINTY LEVEL ($t = 21$, $/MW)

| $\Lambda^\Delta$ | $\Lambda$ | Bus1 | Bus2 | Bus3 | Bus4 | Bus5 | Bus6 |
|---|---|---|---|---|---|---|---|
| 2 | 1   | $-17.666$ | 0 | 0 | 0 | 0 | 0 |
| 2 | 0.8 | $-17.474$ | 0 | 0 | 0 | 0 | 0 |
| 2 | 0.5 | 0 | 0 | 0 | 0 | 0 | 0 |
| 1 | 1   | $-17.474$ | 0 | 0 | 0 | 0 | 0 |
| 1 | 0.8 | $-17.474$ | 0 | 0 | 0 | 0 | 0 |
| 1 | 0.5 | 0 | 0 | 0 | 0 | 0 | 0 |

TABLE XVII
UMP SETTLEMENT W.R.T. UNCERTAINTY LEVEL ($t = 21$, $)

| $\Lambda^\Delta$ | $\Lambda$ | $\Theta_{1,t}^G$ | $\Theta_{2,t}^G$ | $\Theta_{3,t}^G$ | $\Psi_{1,t}$ | $\Psi_{3,t}$ | Revenue |
|---|---|---|---|---|---|---|---|
| 2 | 1   | 780.82 | 178.42 | 60.52 | 1,013.43 | 138.23 | 131.9 |
| 2 | 0.8 | 419.38 | 0      | 12.98 | 435.45   | 11.61  | 14.71 |
| 2 | 0.5 | 0      | 0      | 0     | 0        | 0      | 0 |
| 1 | 1   | 419.38 | 0      | 0     | 544.32   | 14.51  | 139.45 |
| 1 | 0.8 | 419.38 | 0      | 0     | 435.45   | 11.61  | 27.69 |
| 1 | 0.5 | 0      | 0      | 0     | 0        | 0      | 0 |

The optimal generation schedules are presented in Table XIV. With $\Lambda^\Delta = 2$, the generation output of G1 increases from 195.193 to 201.59 MW when $\Lambda$ decreases from 1 to 0.8, while the outputs of G2 and G3 both decrease. As the fuel cost of G1 is lower, the total operation cost is reduced, which is consistent with the data in Table XII. On the other hand, the upward generation reserve (i.e. $Q_1^{\text{up}}$) provided by G1 is reduced by 24-18.411 = 5.599 MW. Similar results can also be observed when the budget parameter $\Lambda^\Delta$ and $\Lambda$ are at other values. Another observation is that if $\Lambda = 0.5$, the generation schedules remain the same when changing $\Lambda^\Delta$. This is also consistent with the total operation cost data shown in Table XII. The generation reserves are the byproducts of the optimal solution to SCED problem.

The upward and downward UMPs are shown in Table XV and Table XVI, respectively. The UMPs are zeros with $\Lambda = 0.5$. It indicates that the small perturbation of the uncertainty at Bus 1 and Bus 3 does not change the total operation cost for the base case. In this case, although both the upward and downward generation reserves are available according to Table XIV, units do not get any credit. In the meantime, the uncertainty sources also do not need to pay. Those generation reserves are "free" byproducts of the SCED solution. As shown in Table XV, by fixing $\Lambda^\Delta = 2$, the UMPs increase with the increment of the uncertainty level $\Lambda$ at Hour 21. For example, the upward UMP is 0 at Bus 1 with $\Lambda = 0.8$. By changing $\Lambda$ to 1, the upward UMP at Bus 1 increases to $14.868/MW. It is also observed that UMPs vary with locations. As shown in Table XVI, the downward UMPs at Bus 1 are non-zero, which indicates that the downward generation reserves at Bus 1 are scarce resources; the downward UMPs at other buses are zeros, which indicates that the downward generation reserves at those buses are not scarce resources. The sensitivity data in these tables show that the uncertainties are only charged when the uncertainty level is above a certain threshold. Accordingly, the flexible resources are also entitled to credits for the uncertainty management.

The market settlement related to uncertainties are shown in Table XVII. The payment to a generator is denoted as $\Theta_{i,t}^G$, and the charge to an uncertainty source is denoted as $\Psi_{m,t}$. They are calculated based on Eqs. (21) and (16) [1], respectively. For example, the payment to G1 is

$$\Theta_{1,21}^G = 14.868 \times 24 + (-17.666) \times (-24) = \$780.82$$

with $\Lambda^\Delta = 2, \Lambda = 1$. $14.868/MW and $-17.666/MW are the upward and downward UMPs, respectively. 24 MW and -24 MW are the upward and downward generation reserves, respectively.

One observation from Table XVII is that G3 gets more reserve credit in the case of $\Lambda^\Delta = 2, \Lambda = 1$ than $\Lambda^\Delta = 2, \Lambda = 0.8$, although it has smaller generation reserves in the first case. The reason is that the upward UMP at Bus 6 soars to $17.493/MWh from $2.596/MWh when $\Lambda$ is increased from 0.8 to 1. The charge to the uncertainty source at Bus 3 is calculated as

$$\Psi_{3,21} = 16.634 \times 8.31 + 0 \times (-8.31) = \$138.23,$$

where $16.634/MW is the upward UMP and 8.31MW is the upper bound of the uncertainty.

Another observation from Table XVII is that the total money uncertainty sources pay decreases with the uncertainty level $\Lambda$. For example, when $\Lambda^\Delta = 2$, the payment of uncertainty source at Bus 1 is $1013.43 at Hour 21 with $\bar{u}_{1,21} = 31.15$MW and $\bar{u}_{3,21} = 8.31$MW. If the uncertainty sources reduce the uncertainties $\bar{u}_{1,21}$ and $\bar{u}_{3,21}$ to $31.15 \times 0.8 = 24.92$MW and $8.31 \times 0.8 = 6.648$MW, respectively, then the payment from uncertainty source at Bus 1 decreases from $577.98 to $435.45, and the uncertainty source at Bus 3 also pay $126.62

TABLE XVIII
G3 WITH MODIFIED CAPACITY AND RAMPING ($t = 21$, $\Lambda^\Delta = 2, \Lambda = 0.8$)

| $P_{3,t}$ | $Q_3^{\text{up}}$ | $Q_3^{\text{down}}$ | $\pi_6^{\text{e}}$ | $\pi_6^{\text{u,up}}$ | $\pi_6^{\text{u,down}}$ |
|---|---|---|---|---|---|
| 15 | 5.5 | -5 | 35.2341 | 2.5961 | 0 |

TABLE XIX
LMP AND UMP WITH MODIFIED ENERGY BID
($\Lambda^\Delta = 2, \Lambda = 0.8, t = 21$, $/MWH)

| | Bus1 | Bus2 | Bus3 | Bus4 | Bus5 | Bus6 |
|---|---|---|---|---|---|---|
| $\pi_{6,21}^{\text{e}}$ | 15.1640 | 30.6380 | 32.1846 | 40.3353 | 38.7887 | 32.9369 |
| $\pi_{6,21}^{\text{u,up}}$ | 0.0000 | 0.0000 | 1.5466 | 9.6973 | 8.1507 | 2.2989 |
| $\pi_{6,21}^{\text{u,down}}$ | -15.4740 | 0.0000 | 0.0000 | 0.0000 | 0.0000 | 0.0000 |

less. It indicates that uncertainty sources have the incentives to reduce the uncertainties.

Next, we check whether the generators have the incentives to increase the ramping rates and generation capacities. Assume that the capacity of G2 is increased to 20.5 MW from 20MW, and the upward ramping is increased to 5.5MW/h from 5MW/h. The optimal schedules and marginal prices are listed in Table XVIII when $\Lambda^\Delta = 2, \Lambda = 0.8$. The new credit for managing uncertainty is

$$5.5 \times 2.5961 = \$14.28,$$

which is larger than the original credit of $12.98. Therefore, G2 is entitled to more generation reserve and credit by increasing its flexibility. However, it also should be pointed out that the market participants are assumed as price takers in this paper [8]–[10]. If the changes of the flexibilities are large enough to change the UMPs, then the market participants may be able to manipulate the market.

In this section, generators only submit the energy bids (i.e. step-wise fuel costs shown in Table IV), and the generation reserve bids are assumed to be zeros. It is also consistent with what is available in the literature [6]. In this paragraph, we demonstrate the impacts of energy bid changes on UMP. Consider the case where $\Lambda^\Delta = 2, \Lambda = 0.8$. The bidding prices of G2 in all segments are reduced by $2/MWh to $30.638/MWh, $30.674/MWh, $30.71/MWh, $30.746/MWh, and $30.782/MWh, respectively, while the bids of other units remain the same. Accordingly, we get another set of price signals after solving the problem (RSCUC). At Hour 21, the generation outputs remain the same. But both LMPs and UMPs are changed as shown in Table XIX. By comparing the LMPs in Table XIX and Table XIII, we observe that the energy prices for all the buses are lower when marginal prices of G2 are reduced. The UMPs also show the same trends by comparing the UMPs in Table XIX, Table XV, and Table XVI. For example, the upward UMP at Bus 4 is reduced to $9.6973/MW from $10.951/MW, the upward UMP at Bus 6 is also reduced by $0.298/MW, and the downward UMP at Bus 1 decreases to $15.474/MW from $17.474/MW. It can be observed that the highest UMP has the largest change in this case. Similar to the discussion in the above paragraph, the trends may be altered if the participants have enough market power. The analysis on Nash-equilibrium would be an interesting future research.

### D. FTR Underfunding

When $\Lambda^\Delta = 2, \Lambda = 1$, the generation schedules at Hour 21 are presented in Table XIV. The power flow of Line 2 is 97.63MW, which is smaller than its physical limit of 100MW. The fact that LMPs are different at different locations indicates that the congestion components in LMPs are non-zero. This is different from the traditional practice that congestion components exist only when certain power flows reach their limits. This is because the transmission reserve 100MW - 97.63MW = 2.47MW is kept to guarantee the delivery of the generation reserve. At least one of the transmission constraints (9f) for Line 2 is binding in problem (RSCED) of the companion paper [1], and at least one of their Lagrangian multipliers is positive. The positive Lagrangian multiplier leads to the congestion component of LMPs.

As stated in Section III-C in [1], the LMPs obtained according to its definition in the robust optimization framework may lead to FTR underfunding. Consider a set of FTR amounts $[202.3429, 23.2771, -55.772, -94.924, -94.924, 20]$. It can be verified that the FTR amounts satisfy the SFT in the FTR market. Then the total credit for the FTR holders is

$$100\text{MW} \times \$55.5477/\text{MW} = \$5,554.77,$$

where $55.5477/MW is the sum of shadow prices for the line constraints, as shown in the proof of Theorem 1 in [1]. However, the congestion fee in the DAM is

$$97.6254\text{MW} \times \$55.5477/\text{MW} = \$5,422.87.$$

It means that the congestion payment collected from DAM is not enough to cover the FTR credit. The FTR underfunding value is

$$\$5554.77 - \$5422.87 = \$131.90.$$

Therefore, the FTR underfunding may occur in this scenario if the existing market mechanism is used. This example shows that the robust solution may lead to FTR underfunding in the traditional market framework.

Now we consider the revenue residue after the UMP settlement in the market. The revenue residues are shown in the last column in Table XVII. When $\Lambda^\Delta = 2, \Lambda = 1$, the revenue residue is $131.9, which is obtained as

$$(1013.43 + 138.23) - (780.82 + 178.42 + 60.52) = \$131.9.$$

It can be observed that the FTR underfunding can be covered by the revenue residue exactly in this scenario. Therefore, the revenue is adequate at Hour 21.

As the UMPs vary with locations, there are congestion components when $\Lambda = 1$ and $\Lambda = 0.5$. According to the data in Table XVII, the revenue residues in these cases are also non-zero. This verifies Eqs. (12) and (23). The fact that the UMPs vary with locations indicates that at least one of the shadow prices $\eta_{l,t}^k$ for line constraints is positive according to (12). The positive $\eta_{l,t}^k$ must lead to the transmission reserve credit from (23). In other words, the revenue residue related to uncertainties must be positive. The revenue residue is then distributed to the FTR holders. Note that **if and only if UMPs vary with locations, the revenue residue related to uncertainties is positive**.

TABLE XX
OPERATION COST AND PAYMENT RELATED TO UMP ($,$\Lambda^{\Delta} = 10$)

| $\Lambda$ | Op. Cost | Un. Payment | Gen. Res. Credit | Rev. Res. |
|---|---|---|---|---|
| 0.2 | 1,866,023 | 11,043 | 10,560 | 483 |
| 0.25 | 1,871,364 | 20,044 | 19,209 | 835 |
| 0.3 | 1,877,471 | 30,879 | 28,658 | 2,221 |

### III. IEEE 118-BUS SYSTEM

The simulations are performed for the IEEE 118-bus system with 54 thermal units and 186 branches in this section. The peak load is 6600MW. The detailed data including generator parameters, line reactance and ratings, and load profiles can be found at http://motor.ece.iit.edu/Data/RSCUC_118_UMP.xls. Two cases are studied in this section.

1) The uncertainty levels and load levels are changed to analyze the simulation results in the system level.
2) An energy storage is installed at a specified bus with high UMP to show the potential application of UMPs.

#### A. Case 1

We assume that the uncertainty sources are located at buses (11, 15, 49, 54, 56, 59, 60, 62, 80, 90). The budget parameter $\Lambda^{\Delta}$ is set to 10 in this section. The bus-level uncertainty budget parameter $\Lambda$ changes from 0.2 to 0.3, and the bound of the uncertainty is the base load. The simulation results are shown in Table XX. It can be observed that the total operation cost increases with increasing $\Lambda$. The largest operation cost is $1,877,471 with $\Lambda = 0.3$. It indicates that a larger uncertainty level may increase the operation cost. The payments and credits related to UMPs are also listed in the table. The columns "Un. Payment" and "Gen. Res. Credit" denote the total payment from uncertainty sources and credit to generation reserves, respectively. It can be observed that the uncertainty sources pay less when the uncertainty level is low. The lowest payment is $11,043 and the highest one is $30,879. On the other hand, the credit entitled to the generation reserves is also a monotonically increasing function of $\Lambda$. When $\Lambda = 0.3$, the generation reserves have the highest credit. The last column "Rev. Res." shows that the revenue residues related to UMPs. It can be observed that the residue is always positive. When $\Lambda$ is large, the residue is also high. As discussed in Section II-D, the residue exists when there are line congestions.

Fig. 3 depicts the heat map for the upward UMPs from Bus 80 to Bus 100 in 24 hours. The x-axis represents time intervals and the y-axis represents bus numbers. The color bar on the right shows different colors for various UMP values. For example, the $0/MW is denoted by the blue color at the bottom, and the $18/MW is represented by the dark red color on the top of the color bar. It can be observed that the uncertainty sources have system-wide unique UMPs at some intervals, such as Hours 8, 13, 15, and so on. It indicates that there is no transmission reserve in these hours. On the other hand, the UMPs at Hour 11 vary dramatically with different locations. The highest upward UMP is around $18/MW, and the lowest one is around $2/MW. According to the data shown in Fig. 3, the high UMP at Bus 94 may attract investment of flexible resources, such as energy storages, in terms of generation reserve credit, and Bus 100 is an attractive location for the investment of renewable energy sources in terms of uncertainty payments.

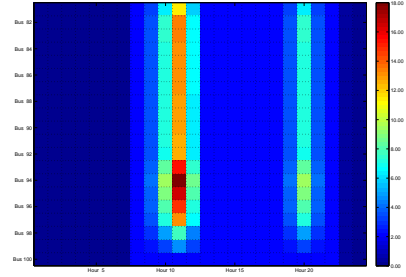

Fig. 3. Heat Map for Upward UMPs ($\Lambda = 0.3$)

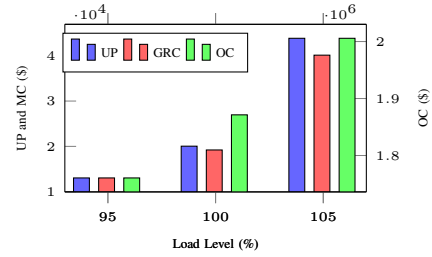

Fig. 4. Uncertainty Payment (UP), Generation Reserve Credit (GRC), and Operation Cost (OC) with Different Load Levels ($\Lambda = 0.25$)

Fig. 4 shows the uncertainty payment and generation reserve credit with respect to load levels. The base load level is set at 100%. Higher loads in general lead to more uncertainty payments and generation reserve credits. It is also consistent with the heat map of UMPs in Fig. 3, where UMPs at peak load hours are high. It suggests that the generation reserves also become scarce resources when load levels are high.

#### B. Case 2

As discussed in Case 1, the upward UMP on Bus 94 is high at Hour 11. Assume that an energy storage (8MW/30MWh) is installed at Bus 94. A simple model for the energy storage is formulated as

$$E_t = E_{t-1} + \rho^d P_t^D + \rho^c P_t^C, \forall t$$
$$0 \leq E_t \leq E^{\max}, \forall t$$
$$0 \leq -P_t^D \leq I_t^D R^D, \forall t$$
$$0 \leq P_t^C \leq I_t^C R^C, \forall t$$
$$I_t^D + I_t^C \leq 1, \forall t$$
$$E_{N_T} = E_0,$$

where $E_t$ denotes the energy level, $P_t^D$ and $P_t^C$ represent the discharging and charging rates, and $I_t^D$ and $I_t^C$ are the indicators of discharging and charging. As the UMP is the major concern in this section, we use simplified parameters for storage. The discharging efficiency $\rho^d$ and charging efficiency $\rho^c$ are set to 100%. The capacity $E^{\max}$ and initial energy level $E_0$ are set to 30 MWh and 15 MWh, respectively. The maximal charging rate $R^D$ and discharging rate $R^C$ are set to 8 MW/h.

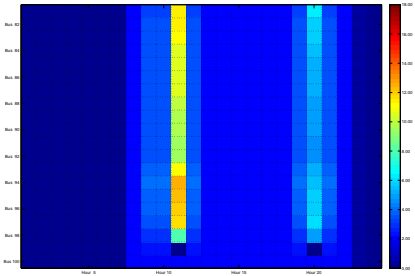

Fig. 5. Heat Map for Upward UMPs with Energy Storage ($\Lambda = 0.3$)

By introducing the energy storage, the new operation cost is reduced to $1,875,211 from $1,877,471. The payment collected from the uncertainty sources becomes $27,473, and the credit to generation reserves decreases to $24,289. Compared to the data in Table XX, the energy storage also helps to reduce the payment related to UMPs. The storage is entitled to $1326 generation reserve credit. Fig. 5 depicts the new upward UMPs after the installation of the energy storage. Compared to that in Fig. 3, the upward UMP for Hour 11 at Bus 94 decreases a lot. The UMPs for Hour 10 and 12 are also lower. It suggests that sitting the energy storage at Bus 94 effectively lower the generation reserve price.

The simulation results show that the energy storage can provide additional upward and downward generation reserves. It is demonstrated that flexible resources can lower the UMPs, and UMPs provide the investment signal at locations where generation reserves are scarce resources.

## IV. CONCLUSION

The two-part paper propose a new market scheme to charge the uncertainty sources and credit the generation reserve providers according to the UMPs in day-ahead electricity markets. The UMP formulation is derived within a robust optimization framework. We also characterize the market equilibrium for the new market mechanism. Our study also shows that traditional pricing mechanism within RSCUC framework may lead to FTR underfunding and the payment collected from uncertainty sources can cover the deficit. Our study shows load serving entities can have lower energy prices within the new market scheme, as the reserve fees are paid by uncertainty sources. It can prevent market participants with uncertainties from gaming the system to certain extent.

An important future research on this topic is that flexible resources, such as generators with large generation reserves and storages, can bid the generation reserves within the RSCUC framework. The Nash-equilibrium for UMP is also our ongoing research. The UMPs derived in this paper also provide an important price signal for the long-term investment of flexible resources. When the upward UMP or downward UMP at a bus is high, the investor can get more return in terms of generation reserves. Many potential applications on UMPs are open for future research.

**Hongxing Ye** (S'14) received his B.S. degree and M.S. degree from Xi'an Jiaotong University, China in 2007 and 2011, both in electrical engineering. He is currently working toward the Ph.D. degree at the Illinois Institute of Technology, Chicago. His research interests include optimization, economic operation and security analysis of power system.

**Yinyin Ge** (S'14) received her B.S. degree and M.S. degree from Xi'an Jiaotong University, China in 2008 and 2011, both in electrical engineering. She is currently a Ph.D. candidate of Electrical Engineering at the Illinois Institute of Technology, Chicago. Her research interests are power system optimization and modeling, Smart Grid and power system stability and control.

**Mohammad Shahidehpour** (F'01) received his Ph.D. degree from the University of Missouri in 1981 in electrical engineering. He is currently the Bodine Chair Professor and Director of the Robert W. Galvin Center for Electricity Innovation at the Illinois Institute of Technology, Chicago.

**Zuyi Li** (SM'09) received the B.S. degree from Shanghai Jiaotong University, Shanghai, China, in 1995, the M.S. degree from Tsinghua University, Beijing, China, in 1998, and the Ph.D. degree from the Illinois Institute of Technology (IIT), Chicago, in 2002, all in electrical engineering. Presently, he is a Professor in the Electrical and Computer Engineering Department at IIT. His research interests include economic and secure operation of electric power systems, cyber security in smart grid, renewable energy integration, electric demand management of data centers, and power system protection.